\documentclass[11pt]{amsart}

\usepackage{amsmath,amsfonts,amsthm,amssymb,amscd,verbatim,graphicx,color,upgreek}
\usepackage[notcite,notref]{showkeys}
\usepackage[margin=2cm]{geometry}

\usepackage{hyperref}

\DeclareMathSymbol{\Gamma}{\mathalpha}{operators}{0}

\newcommand{\SL}{\mathrm{SL}}
\newcommand{\AGL}{\mathrm{AGL}}

\newcommand{\AGammaL}{\mathrm{A\Upgamma L}}

\newcommand{\D}{\mathrm{\mathbf{D}}}
\newcommand{\De}{\mathcal{D}}

\newcommand{\Sym}{\operatorname{Sym}}

\newcommand{\supp}{\operatorname{supp}}
\newcommand{\Alt}{\operatorname{Alt}}

\newcommand{\PSL}{\operatorname{PSL}}

\renewcommand{\Gamma}{\varGamma}
\renewcommand{\epsilon}{\varepsilon}

\renewcommand{\leq}{\leqslant}
\renewcommand{\geq}{\geqslant}

\newcommand{\B}{\mathcal{B} }

\newtheorem*{thma}{Theorem A}
\newtheorem*{thmb}{Theorem B}
\newtheorem*{thmc}{Theorem C}

\newtheorem{prop}{Proposition}[section]
\newtheorem{thm}[prop]{Theorem}

\newtheorem{lem}[prop]{Lemma}

\theoremstyle{definition}

\numberwithin{equation}{section}

\def\doubleunderline#1{\underline{\underline{#1}}}

\begin{document}

\title{Groups Obtained from $2-(n,4,3)$ Supersimple Designs}

\author{Nick Gill}\thanks{}
\address{Department of Mathematics,
University of South Wales,
Treforest, CF37 1DL
}
\email{nick.gill@southwales.ac.uk}

\author{Jerem\'ias Ram\'irez}
\address{Escuela de Matem\'aticas, 
Universidad Nacional, 
Heredia, 
Costa Rica.
}
\email{jeremias.ramirez.jimenez@una.ac.cr}

\begin{abstract}
We contribute towards the classification programme for Conway groupoids associated to a 
$2-(n,4,\lambda)$ design. Our main results improve the known bounds for a hole stabilizer to be 
primitive, or to contain the alternating group, $\Alt(n-1)$. We exploit these improved bounds to 
give a partial classification for Conway groupoids when $\lambda=3$.
\end{abstract}

\maketitle

\section{Introduction}

In his famous paper \cite{Co1}, John Conway used a ``game'' played on  the projective plane 
${\mathbb P}_3$ of order $3$  to construct the sporadic Mathieu group $M_{12}$, as well as a 
special subset of $\Sym(13)$ which he called $M_{13}$, and which could be endowed with the 
structure of a groupoid.

In recent work (\cite{GGNS, GGS}), Conway's construction has been generalized to geometries other 
than ${\mathbb P}_3$, namely to supersimple $2-(n,4,\lambda)$ designs. In this more general context, 
the analogue of the group $M_{12}$ is a subgroup of $\Sym(n-1)$ which is known as the 
{\it hole stabilizer} of the design. In this paper we prove a number of results concerning 
hole stabilizers. Our first main result is the following. It is a strengthening of \cite[Theorem~E]{GGS}.

\begin{thma}
Suppose that $\De$ is a supersimple $2-(n,4,\lambda)$ design, and that $\infty$ is a point in $\De$. 
Let $G:=\pi_\infty(\De)$ be the hole stabilizer of $\infty$, considered as a permutation group via 
its natural embedding in $\Sym(n-1)$.
\begin{enumerate}
\item If $n>\frac{24}{7}\lambda+1$, then $G$ is transitive;
\item if $n> 9\lambda-6$, then $G$ is primitive;
\item if $n\geq 10\lambda-5$ then $G$ is generously transitive;  
\item  if $n>9\lambda^2-12\lambda+5$, then one of the following holds:
\begin{enumerate}
\item $G$ contains $\Alt(n-1)$;
\item $\lambda=1$, $\De=\mathbb{P}_3$ (the projective plane of order $3$), and $G =  M_{12}$. 
\end{enumerate}
\end{enumerate}
\end{thma}

The value $\frac{24}{7}\lambda+1$ in the first bound in Theorem~A is an improvement on the $4\lambda+1$ 
which appears in \cite[Theorem~E]{GGS}. In fact Lemmas~\ref{l: intransitive} and \ref{l: 2 orbs}, 
which are stated and proved in \S\ref{s: trans}, give even stronger bounds, and the first bound in 
Theorem~A follows directly from these results. The second bound in Theorem~A is proved in 
\S\ref{s: bound prim}.

The third bound in Theorem~A is of a different flavour to results in \cite{GGS}. A definition is 
required: a permutation group $G\leq \Sym(n)$ is called \emph{generously transitive} if for each 
$i,j \in \{1,2,\ldots , n\}$, $i\neq j$ exists an element $g\in G$ that interchanges this elements.

The fourth bound in Theorem~A is already known and appears in \cite[Theorem~E]{GGS}. We keep it as 
part of our Theorem~A as it will be useful later. However we can also give another result which 
has a similar flavour.

\begin{thmb}
Let $G=\pi_{\infty}(\De)$. If $n>18\lambda-17$, then one of the following holds:
\begin{enumerate}
\item $\Alt(\Omega\setminus \{\infty\}) \leq G$; 
\item $\De$ is the projective plane of order $3$ and $G\cong M_{12}$;
\item $G\cong \Sym(m)$ with $m\leq 3\lambda-1$, and the action of $G$ on $\Omega\setminus\{\infty\}$ 
is permutation isomorphic to the action on the set of $k$-subsets of $\{1,\dots, m\}$ for some 
$k\in\{2,\dots, \lfloor\frac{m}{2}\rfloor\}$;
\item $G\cong \Alt(m)$ with $m\leq 2\lambda$, and the action of $G$ on $\Omega\setminus\{\infty\}$ 
is permutation isomorphic to the action on the set of $k$-subsets of $\{1,\dots, m\}$ for some 
$k\in\{2,\dots, \lfloor\frac{m}{2}\rfloor\}$. 
\end{enumerate}
\end{thmb}

The strength of Theorem~B is that it yields conclusions (1) and (2), conditional only a {\bf linear} 
lower bound in $\lambda$ for $n$, as opposed to the quadratic lower bound in Theorem A (4).

On the other hand, the weakness of Theorem~B is in conclusions (3) and (4): these parts have the 
advantage that they explicitly describe the permutation group $G$, however there is an associated 
loss of control on the size of $n$ in terms of $\lambda$. Nonetheless, for many values of $m$ and 
$k$ these actions violate the original quadratic bound that is the third bound in Theorem~A and 
so, in principle, one could use this to obtain further restrictions on $m$ and $k$ in terms of 
$\lambda$. We have elected not to do this as it would introduce ``clutter'' to the statement, 
but practical applications of this theorem will probably require such an analysis. We prove 
Theorem~B in \S\ref{s: b}.

Our final theorem extends the classification of hole stabilizers for small values of $\lambda$. 
When $\lambda=1$ or $2$, \cite[Theorem C]{GGNS} gives a full classification. We now partially 
deal with the case of $\lambda=3$.

\begin{thmc}
 Let $\De=(\Omega, \B)$ be a $2-(n,4,3)$ supersimple design. Let $\infty\in\Omega$ 
and set $G:=\pi_\infty(\De)$. Then either $G$ contains $\Alt(\Omega\setminus\{\infty\})$ or one of 
the following holds:
\begin{enumerate}
\item $G$ is {\bf intransitive}: $(n,G)=\doubleunderline{(8, \{1\})}$;
\item $G$ is {\bf imprimitive}: $(n,G)$ is in
\[
\{ \doubleunderline{(9, \Alt(4) \wr C_2)}, (13, ?), (16, ?), (17, ?), (21, ?)\}
\]
\item $G$ is {\bf primitive} $(n,G)$ is in
\[
\left\{ \begin{array}{l} (13, M_{12}), (13, M_{11}), (16, \SL_4(2)), \underline{(16, \Sym(6))}, \\ (16, \Alt(7)), (16, \Alt(6)), (17, ?)\end{array} \right\}.
\]
\end{enumerate}
Furthermore, for those entries that are double underlined, all such examples are known and classified; for those entries that are single underlined, an example is known; for those entries that are not underlined, no such example is known.
\end{thmc}

Note that Theorem~C asserts that the classification of hole stabilizers is complete for $\lambda=3$ except when $n\in \{13,16,17, 21\}$. In \S\ref{s: l3} we prove Theorem~C, and give a full description of all relevant examples.

One final remark: as we have said, our three main results fit into the programme of classification for Conway groupoids. In fact, though, we never directly study the groupoid of a design itself -- all of our results are stated in terms of ``the hole stabilizer'', and our proofs are also couched in these terms. For a definition of the Conway groupoid associated with a supersimple $2-(n,4,\lambda)$ design we refer to \cite{GGS} where, in addition, the connection between the hole stabilizer and the Conway groupoid is made clear.

\subsection{Acknowledgments} It is a pleasure to thank B.~McKay and M.~Meringer who very kindly did a number of computer calculations at our request.

\section{Background}

\subsection{Block designs}\label{s: block design}

Let $t,n,k,\lambda$ positive integers. A \emph{balanced incomplete block design} $(\Omega, \mathcal{B})$, 
also known as a \emph{$t-(n,k,\lambda)$ design}, is a finite set $\Omega$ of size $n$, together 
with a finite multiset $\mathcal{B}$ each of size $k$ (called \emph{lines}) such that any subset 
of $\Omega$ of size $t$ is contained in exactly $\lambda$ lines.

In this paper we are mostly interested in $2-(n,4,\lambda)$ designs. Such a design is called 
\emph{simple} if there are no repeated lines, and \emph{supersimple} if any two lines intersect 
in at most two points. In what follows we will be interested exclusively in supersimple 
$2-(n,4,\lambda)$ designs and so we can assume that the multiset $\mathcal{B}$ is in fact a set.

For some values of $t,n,k,\lambda$ the set of $t-(n,k,\lambda)$ designs have been completely enumerated.
 We will use this information for computer calculation purposes. We refer the reader to \cite{Col} 
for more information.

Let us note a particularly important example: the \emph{Boolean quadruple system of order $2^k$} is 
the design $(\Omega, \mathcal{B})$, where $\Omega$ is identified with the set of vectors of $\mathbb{F}_2^{k}$, 
and 
\[\mathcal{B}:=\{ \{ v_1,v_2,v_3,v_4 \} :v_i \in \Omega \quad \mathrm{and} \quad v_1+v_2+v_3+v_4=0 \}.\]

It is easy to see that $\De$ is a $2-(2^k,4,2^{k-1}-1)$ design (in particular, when $k=2$, 
it is a $2-(8,4,3)$ design). 


\subsection{Permutation Groups}\label{s: permutation groups}

In this subsection, we collect some related notions about permutation groups that will be used at 
long of this paper. For more details we refer the reader to \cite{Mor}.

Suppose that $G$ is a group acting on a non-empty set $\Omega$. The action is called \emph{transitive} 
if for all $x,y \in \Omega$ there is an element $g\in G$ such that $x^g=y$. 

Suppose that the action of $G$ on $\Omega$ is transitive. A \emph{system of imprimitivity} is a 
partition of $\Omega$ into $l$ subsets $\Delta_1, \Delta_2, \ldots , \Delta_l$ each of size $k$ such 
that $1<k,l<n$, and so that for all $i\in \{1,2,\ldots , l\}$ and all $g\in G$ there exists 
$j\in \{1,2,\ldots l\}$ such that 

\[\Delta_i^g=\Delta_j.\]

The sets $\Delta_i$ are called \emph{blocks}. We say that G acts \emph{imprimitively} if there 
exists a system of imprimitivity. If no such set exists then G acts \emph{primitively} on $\Omega$.

The \emph{support} of an element $g\in G$, denoted $\supp(g)$ is the set of points in 
$\Omega$ not fixed by $g$.

\subsection{Hole Stabilizers}

Suppose that $\De=(\Omega, \B)$ is a supersimple $2-(n,4,\lambda)$ design. Two points $x,y \in \Omega$ 
are \emph{collinear} if there is some line in $\mathcal{B}$ that contains $x$ and $y$.

Suppose that a pair of distinct elements $x,y\in \Omega$ are collinear. We define the \emph{elementary move} 
$[x,y]$ to be the permutation

\[[x,y]:=(x,y)\prod_{i=1}^{\lambda}(x_i,y_i),\]

where $\{x,y,x_i,y_i\}$ is a line in $\mathcal{B}$ for every $1\leq i\leq \lambda$. This product is 
well defined because $\De$ is supersimple. We also define $[x,x]=\mathrm{Id}_{\Omega}$.

Let $a$ and $b$ be distinct points in $\Omega$. We define

\[
\overline{a,b}:=\{x\in\Omega\,|\,\textnormal{there exists $\ell\in\B$ such that $x,a,b,\in\ell$}\}.
\]

In particular, note that $a,b\in\overline{a,b}$. Clearly the set of points in $\Omega$ moved by the 
permutation $[a,b]$ (also called the \textit{support} of $[a,b]$) is precisely the set $\overline{a,b}$.

A \emph{move sequence} is

\[[a_0,a_1,a_2, \ldots , a_n]=[a_0,a_1][a_1,a_2][a_2,a_3]\ldots[a_{n-1},a_n]\]

where $a_{i},a_{i+1}$ are collinear for all $0\leq i \leq n-1$. A move sequence is called \emph{closed} 
if $a_0=a_n$. For each $x\in \Omega$ we define the \emph{hole stabilizer}, $\pi_x(\De)$, to be set of all 
closed move sequences such that $a_0=a_n=x$, that is

\[\pi_{x}(\mathcal{D}):=\{[a_0,a_1,\ldots ,a_n]:a_0=a_n=x\}.
\]

It is easy to check that $\pi_{x}(\mathcal{D})$ is a subgroup of $\Sym(\Omega\setminus\{x\})=\Sym(n-1)$. 
In what follows we will need two easy facts \cite[Lemma 3.1 and Theorem A]{GGNS}.

\begin{lem}\label{Properties of D}
Suppose that $\De=(\Omega, \B)$ is a supersimple $2-(n,4,\lambda)$ design and that $x,y\in \Omega$.
\begin{enumerate}
\item $\pi_{x}(\mathcal{D})=\langle [x,a,b,x]:a,b\in \Omega\setminus \{x\} \rangle$.
\item $\pi_x(\mathcal{D})$ and $\pi_y(\mathcal{D})$ are conjugate subgroups of $\Sym(\Omega)$.
\end{enumerate}
\end{lem}

The second statement above implies that all hole stabilizers for a supersimple design $\De$ are 
permutation isomorphic groups. This allows us to talk of ``the'' hole stabilizer $\D$ (defined up 
to permutation isomorphism), and in the rest of this paper we denote this group as $\pi_{\infty}(\mathcal{D})$.

\section{A proof of Theorem~A}

In this section we prove Theorem~A. Throughout this section $\De$ is a supersimple $2-(n,4,\lambda)$ 
design with point set $\Omega$ and $\infty$ one such point. We write $G=\pi_\infty(\De)$. 

\subsection{A bound for transitivity}\label{s: trans}

The lemmas in this section immediately yield statement (1) in Theorem~A.

\begin{lem}\label{l: intransitive}
Suppose that $G=\pi_\infty(\De)$ has $t$ orbits on $\Omega\setminus\{\infty\}$ with $t>1$. Then
\[
 n\leq \frac{2t\lambda}{t-1}+1.
\]
In particular, if $n>4\lambda+1$, then $G=\pi_\infty(\De)$ is transitive.
\end{lem}
\begin{proof}
 Suppose that $x\in \Omega\setminus\{\infty\}$, and write $\Delta_x$ for the orbit of $x$ under $G$. Observe that
 \begin{equation}\label{e: b}
   \Delta_x \supseteq (\overline{\infty, x})^c.
 \end{equation}
Choose $x$ so that $\Delta_x$ is as small as possible. Then $|\Delta_x|\leq \frac{n-1}{t}$. Now, by the observation above,
\[
 \Omega\setminus\{\infty\} = \left(\Delta_x \cup \overline{\infty, x}\right)\setminus\{\infty\}.
\]
Noting that $x\in \Delta_x\cap \overline{\infty, x}$, we obtain that
\[
 n-1 \leq \frac{n-1}{t} + 2\lambda + 1 - 1.
\]
Rearranging this inequality gives the result. The ``in particular'' part of the lemma follows by taking $t=2$.
\end{proof}

\begin{lem}\label{l: 2 orbs}
Suppose that $G=\pi_\infty(\De)$ has $2$ orbits on $\Omega\setminus\{\infty\}$. Then $n\leq \frac{24}{7}\lambda+1$.
\end{lem}

\begin{proof}
Note that if $\lambda\leq 2$, then \cite[Theorem~C]{GGNS} implies that $G=\pi_\infty(\De)$ is always transitive. Thus we assume that $\lambda \geq 3$.

Suppose that $x\in \Omega\setminus\{\infty\}$. Write $\Delta_x$ for the orbit of $x$ under $G$, and note that \eqref{e: b} still holds. Suppose that  $y\in \Omega\setminus(\Delta_x \cup\{\infty\})$. Now, taking complements of both sides of \eqref{e: b} we observe that
\[
\Delta_y=\Delta_x^c \subseteq \overline{\infty, x}.
\]
Note, too, that \eqref{e: b} implies that $\Delta_x\cup \overline{\infty, x}=\Omega$. We wish to give a lower bound for $\Delta_x\cap\overline{\infty, x}$. 

 Observe that there are $\frac{\lambda(n-1)}{3}$ lines through  $\infty$. All of these lines have either at least two elements of $\Delta_x$ or at least two of $\Delta_y$. 
Choose $x$ so that at least half of them (i.e.\ at least $\frac{\lambda(n-1)}{6}$ of them) contain at least 
two elements of $\Delta_x$.

Define
\[
\Lambda=\{(x_1,y) \mid x_1,y\in \Delta_x, \, x_1\neq y, \, y\in \overline{x_1,\infty} \}.
\]
Counting this in two different ways, we obtain that
\[
\Delta_x \cdot (\textrm{Average number of points in $\Delta_x \cap (\overline{x_1,\infty}\setminus\{x_1\})$}) \geq \frac{\lambda(n-1)}{6} \cdot 2
\]
and so we conclude that there exists an element $x$ such that
\[ | \Delta_x \cap \overline{x,\infty} | \geq \frac{\lambda(n-1)}{3\Delta}+1,\]
where $\Delta = |\Delta_x|$. This means that
\[
 n = |\Delta_x \cup \overline{\infty, x}| \leq \Delta+ 2\lambda+2 - \left(\frac{\lambda(n-1)}{3\Delta}+1\right).
\]
Rearranging we obtain that
\[
 n\leq \Delta+1 + \frac{5\Delta\lambda}{3\Delta+\lambda}.
\]
Now, for fixed $\Lambda$, the function $\frac{5\Delta\lambda}{3\Delta+\lambda}$ is an increasing function in the variable $\Delta$. Since $\Delta_x\subset \overline{\infty, y}$, we know that $\Delta\leq 2\lambda$ and we obtain that then
\[
 \frac{5\Delta\lambda}{3\Delta+\lambda}\leq \frac{10}{7}\lambda,
\]
and we obtain that $n\leq \Delta+1+\frac{10}{7}\lambda \leq \frac{24}{7}\lambda+1$.
\end{proof}


\subsection{A bound for primitivity}\label{s: bound prim}

In this section we prove statement (2) of Theorem~A.
Throughout this section we suppose that $G$ is transitive and preserves a system of imprimitivity with $\ell$ blocks each of size $k$ (so that $n-1=k\ell$). Let us start with the following lemma which is \cite[Lemma 6.2]{GGNS}.

\begin{lem}\label{l: prim}
If $n > 4\lambda+1$, then at least one of the following holds:
\begin{enumerate}
\item[(i)] if $a_1,a_2\in \Omega$ lie in the same block of imprimitivity, then $\infty\in\overline{a_1,a_2}$;
\item[(ii)] $n\leq \frac{6\ell}{\ell-1}\lambda+1$.
\end{enumerate}
\end{lem}

Let us label blocks in the system of imprimitivity by $A,B,C,\dots$. Now we label points in $A$ by 
$a_1,a_2,a_3,\dots$, points in $B$ by $b_1,b_2,b_3,\dots,$ and so on.

\begin{lem}\label{l: prim 2}
Suppose that there exists a line $\{a_1,a_2,b,\infty\}$. Then $n\leq \frac{\ell}{\ell-1}(6\lambda-7-\frac{1}{\ell})$.
\end{lem}
\begin{proof}
Choose $x$, a point in $\Omega$ such that $x\not\in \overline{\infty, a_1}\cup \overline{a_1, b} \cup \overline{\infty, b}$. Let $g_x=[\infty, a_1, x, \infty]$ and observe, first, that $a_2^{g_x}=b$. Thus $A^{g_x}=B$. 
Observe, second, that $a_1^{g_x}=x$ and so $x\in B$. We conclude that
\[
B\supseteq \Omega\setminus \left(\overline{\infty, a_1}\cup \overline{a_1, b} \cup \overline{\infty, b}\right) \cup\{b\}.
\]
In particular, $|B|\geq n-(6\lambda-6)+1$. Now use the fact that $|B|=\frac{n-1}{\ell}$, 
and the result follows.
\end{proof}

\begin{lem}\label{l: ell 2}
Suppose that $G$ preserves a system of imprimitivity with $\ell=2$ blocks of size $\frac{n-1}2$. 
Then $n\leq 6\lambda+3$.
\end{lem} 
\begin{proof}
This implies that $G$ contains an element of support of size $2k=n-1$ in any generating set. Now 
the result follows from the fact that $G$ is generated by elements with support of size at most 
$6\lambda+2$ (\cite[Lemma 7.3]{GGNS} -- or see item (4) of Lemma~\ref{Properties of D}).

\end{proof}

The following lemma is stated for $\ell=3$; it is possible that similar statements may hold more generally.

\begin{lem}\label{l: ell 3}
 Suppose that $\ell=3$ and that any line containing $\infty$ contains points from all blocks of imprimitivity ($A$, $B$ and $C$). Then $n\leq 9\lambda-8$.
\end{lem}
\begin{proof}
Let $L=\{\infty, a, b,c\}$ be a line. Then observe that
\begin{align*}
[\infty, a]=(b,c)(b_1, c_1)(b_2, c_2)\cdots (b_{\lambda-1}, c_{\lambda-1}); \\
[\infty, b]=(a,c)(a_1, c'_1)(a_2, c'_2)\cdots (a_{\lambda-1}, c'_{\lambda-1}); \\
\end{align*}
Now consider the element $g=[\infty, a, b,\infty]$. If $g$ is to fix $B$ set-wise, then 
$[a,b]$ must move $c, c_1,\dots, c_{\lambda-1}$ to elements in $B$. If this is the case, then $g$ must interchange $A$ and $C$. The same argument works if we consider what happen when we fix $A$ or $C$ set-wise. 

We conclude that in any case $g$, which is an element of support at most $6\lambda-6$, must move 
at least two blocks, and so
\[
\frac23(n-1)\leq 6\lambda-6.
\]
\end{proof}

\begin{lem}\label{l: ell 4}
Suppose that $L$ is any line containing $\infty$, then $L$ intersects a block of 
imprimitivity in at most $1$ point. Then a block of imprimitivity has size at most 
$2\lambda-1$.
\end{lem}
\begin{proof}
Let $L=\{\infty, a, b,c\}$ and suppose that $x\in B$, $x\not\in\overline{a,\infty}$ 
and $x\not\in\overline{a,b}$. Observe that the supposition implies that 
$x\not\in \overline{b,\infty}$. Now let $g_x=[\infty, a,x,\infty]$ and observe 
that $a^{g_x}=x\in B$ and $c^{g_x}=b\in B$. This is a contradiction.
 
Thus either $x\in\overline{a,\infty}$ or $x\in\overline{a,b}$. The supposition 
ensures that $|B\cap \overline{a,\infty}|\leq \lambda$. Suppose, then that 
$x\in\overline{a,b}$ and $x\not\in\overline{a,\infty}$. Then there is a line 
$\{a,b,x, y\}$ and, defining $g_x$ as before, observe that $a^{g_x}=x\in B$ and 
$c^{g_x}=y^{[x,\infty]}$. If $y\in B$, then the supposition guarantees that 
$c^{g_x}=y^{[x,\infty]}=y\in B$, which is a contradiction. We conclude that 
$y\not\in B$. Thus $\overline{a,b}$ can contain at most $\lambda-1$ points of 
$B$ apart from $b$. The result follows. 
\end{proof}

Let us sum up the work of this section with the next lemma which is statement 
(2) of Theorem~A.

\begin{lem}\label{l: prim final}
If $G$ preserves a non-trivial system of imprimitivity, then $n\leq 9\lambda-6$.
\end{lem}
\begin{proof}
If $\lambda \leq 2$, then the result follows immediately from \cite[Theorem~C]{GGNS}. Assume from here on that $\lambda \geq 3$. If $\ell=2$, then Lemma~\ref{l: ell 2} implies that $n\leq 6\lambda+3$ and the result follows.

Suppose from here on that $\ell\geq 3$. If there exists a line 
$\{a_1,a_2,b,\infty\}$, then Lemma~\ref{l: prim 2} implies that
\[
 n\leq \frac{\ell}{\ell-1}\left(6\lambda-7-\frac{1}{\ell}\right)\leq 9\lambda-11,
\]
and the result follows.

Suppose from here on that if $L$ is any line containing $\infty$, then $L$ 
intersects a block of imprimitivity in at most $1$ point. If $\ell=3$, then 
Lemma~\ref{l: ell 3} implies that $n\leq 9\lambda-8$, and the result follows. 
If $\ell=4$, then Lemma~\ref{l: ell 4} implies that $n\leq 8\lambda-3$, and 
the result follows.

Suppose from here on that $\ell\geq 5$. Then Lemma~\ref{l: prim} implies that
\begin{equation}\label{ee}
n\leq \frac{6\ell}{\ell-1}\lambda+1 \leq 7.5\lambda+1,
\end{equation}
and the result follows for $\lambda \geq 4$. For $\lambda=3$, \eqref{ee} implies 
that $n\leq 23$. Since $2-(n,4,3)$ designs only occur for $n\equiv 0,1\pmod 4$, 
we conclude that $n\leq  21$, and the result follows.
\end{proof} 

\subsection{A bound for generous transitivity}

In this section we prove the third bound in Theorem~A.

\begin{lem}\label{gt:bound}
Let $G=\pi_{\infty}(\mathcal{D})$. If $n\geq 10\lambda-5$ then $G$ is generously transitive.  
\end{lem}

\begin{proof}
Let $a,b\in \Omega\setminus\{\infty\}$. We must find $g\in G$ such that $a^g=b$. If $\infty \not \in \overline{a,b}$ then we can use $g:=\left[\infty,a,b,\infty\right]$.

Suppose that $\infty \in \overline{a,b}$. This means that there exists a line $\{\infty, a, b, c\}$ for some $c$. Choose $x$ such that $x$ is not in $\overline{a,b}$, $\overline{\infty,a}$, $\overline{\infty,b}$,
$\overline{a,c}$, $\overline{b,c}$. Then we can take $g$ to be $\left[\infty,c,x,\infty\right]$. Finally, observe that the sets listed above together contain at most $10\lambda-6$ elements,
so, assuming $n\geq 10\lambda-5$ we obtain the result.

\end{proof}

\section{A proof of Theorem~B}\label{s: b}

Our aim in this section is to prove Theorem~B. We need two background results. The first 
is \cite[Theorem D]{GGS}.

\begin{thm}\label{Trivial E}
Suppose that $\De$ is a supersimple $2-(n,4,\lambda)$ design, and that $[\infty, a, b, \infty]=(\,)$ 
whenever $\infty$ is collinear with $a,b$. Then one of the following is true:
\begin{enumerate}
\item $\De$ is a Boolean quadruple system and $\pi_\infty(\De)$ is trivial;
\item $\De$ is the projective plane of order $3$ and $\pi_\infty(\De) \cong M_{12}$; or 
\item $\pi_\infty(\De)\supseteq \Alt(\Omega\setminus\{\infty\})$.
\end{enumerate}
\end{thm}

We also need a result of Liebeck and Saxl \cite[Theorem 2]{LS}. 

\begin{thm}\label{t: ls}
Let $G$ be a primitive group of degree $d$. Then either
\begin{enumerate}
\item all non-trivial elements have support at least $\frac13 d$ points, or
\item $G$ is a subgroup of $\Sym(m)\wr \Sym(r)$ containing $(\Alt(m))^r$, with $m\geq 5$, 
where the action of $\Sym(m)$ is on $k$-element subsets of $\{1,\dots, \ell\}$ and the wreath 
product has the product action of degree $d=\binom{\ell}{k}^r$.
 \end{enumerate}
\end{thm}

We note that there is an improvement on Liebeck and Saxl's result due to Guralnick and 
Magaard \cite{GM} -- we have elected not to use their result as it includes a longer list of exceptions.

Recall that the product action of $\Sym(\ell)\wr \Sym(r)$ can be thought of as an action 
on the set of functions $\Delta\to\Gamma$, where $\Delta$ is a set of size $r$ and 
$\Gamma$ is a set of size $\ell$. Let $bg=(b_1,\dots, b_r)g$ be an element of 
$\Sym(\ell)\wr \Sym(r)$ (so $b_1,\dots, b_r\in\Sym(\ell)$ and $g\in \Sym(r)$), 
then for $\alpha:\Delta\to\Gamma$, we have
\[
 \alpha^{(b,g)}: \Delta \to \Gamma, \,\, i \mapsto ( i^{g^{-1}} \alpha )^{b_{i^{g^{-1}}}}.
\]
Note that there are $d=\ell^r$ functions $\Delta\to\Gamma$. Using the notation just 
established, we have the following lemma.

\begin{lem}\label{l: product fixed}
\hspace{2em}
\begin{enumerate}
\item Let $G=\Sym(\ell)\wr \Sym(r)$, considered as a permutation group via the product 
action on $d=\ell^r$ points. Suppose that $g=b h \in G$, with $b \in \Sym(\ell)^r$, $h \in \Sym(r)$ 
and $h\neq 1$. Then the number of fixed points of $b h$ is maximal when $h$ is a 
transposition, and $b=1$. In this case $b h$ fixes $\ell^{r-1}=d/\ell$ points.
\item Let $G=\Sym(m)$ acting on the set $\Lambda$ of $k$-subsets of $\{1,\dots, m\}$ 
for some $k\in\{2,\dots, \lfloor\frac{m}{2}\rfloor\}$. 
\begin{enumerate}
\item If $g\in \Sym(m)\setminus\{1\}$, then $g$ moves at least $2m-4$ points of $\Lambda$.
\item If $g\in \Alt(m)\setminus\{1\}$, then $g$ moves at least $3m-6$ points of $\Lambda$.
\end{enumerate}
\end{enumerate}
\end{lem}
\begin{proof}
For (1), let $h$ be non-trivial, and label elements so that $1^h=2$. If $g=b h$ fixes a function $\alpha$, then we require that
\[
 (1 \alpha) = (2 \alpha) ^{b_2}.
\]
Thus the image of $1$ under $\alpha$ is prescribed by the image of $2$, and we obtain immediately that there are at most $\ell^{r-1}$ possibilities for $\alpha$.

For (2), let $g$ be non-trivial, and label elements so that $1^g=2$. The number of $k$-sets that contain $1$ but don't contain $2$ is $\binom{m-2}{k-1}$; likewise the number of $k$-sets that contain $2$ but don't contain $2^g$ is  $\binom{m-2}{k-1}$. These two families of sets are disjoint, and all sets contained therein are moved by $g$, hence $2\binom{m-2}{k-1}$ is a lower bound on the number of points moved by $g$.

If $k>2$ and $m>5$, then this immediately yields the lower bound $3m-6$ (recall that we may assume that $k\leq m/2$). Thus we must consider the cases $m\leq 5$ or $k=2$; note, though, that if $m\leq 5$, then we automatically have that $k\leq 2$.

Suppose, then, that $k=2$. If in the cycle decomposition of $g$, we have $(1,2,\dots, t)$, then the number of $2$-sets containing $1$ but not $2$, then $2$ but not $3$ (and so on ) is at least $t(m-2)$. Thus if $g$ contains a cycle of length $3$ or more, then the result follows; the bound for (a) also follows. Suppose, then that $g$ is in $\Alt(m)$ and $g$ is a product of $k$ distinct transpositions with $k\geq 2$; write $g=(1,2)(3,4)\cdots$. Then the same argument yields a lower bound of $4 m-8$, and the result follows.
\end{proof}

We are ready to prove Theorem~B.

\begin{proof}[Proof of Theorem~B]
The result is true for $\lambda\leq 2$ by classification theorems in \cite{GGNS}. Note that $n>18\lambda-17>9\lambda+1$ for $\lambda \geq 3$ and so $G=\pi_\infty(\De)$ is a primitive subgroup of $\Sym(n-1)$.

Now, by Theorem~\ref{Trivial E}, we can assume that $[\infty,a,b,\infty]\neq 1$ for some $a,b$ collinear with $\infty$. Such an element has support at most $6\lambda-6$. Now we consider the possibilities given in Theorem~\ref{t: ls}. If possibility (1) occurs, then we conclude that $n=d-1$ with
\[
 \frac13d \leq 6\lambda-6
\]
which is a contradiction.

Thus, possibility (2) occurs: $G$ is a subgroup of $\Sym(\ell)\wr\Sym(r)$ in the product action on $\ell^r$ points. If $r\neq 1$, then any set of generators for $G$ must include an element $bg$ with $g\neq 1$ (using the notation established before Lemma~\ref{l: product fixed}). However we know that the set of elements of the form $[\infty, a,b,\infty]$ generate $\pi_\infty(\De)$ and these elements have support at most $6\lambda+2$. Referring to Lemma~\ref{l: product fixed}, we conclude that $n=d-1$ with
\[
 d-d/\ell\leq 6\lambda+2
\]
and so $d<12\lambda+4$ which is a contradiction for $\lambda \geq 4$. For $\lambda=3$, we have a contradiction when $\ell\neq 2$. When $\ell=2$ we must rule out $n\in \{37, 38, 39,40\}$ but, since none of these are powers of $2$, this is immediate.

Thus we are left with the possibility that $r=1$, $d=\binom{m}{k}$, and $G=\pi_\infty(\De)$ is either $\Sym(m)$ or $\Alt(m)$ with the action on $\Omega\setminus\{\infty\}$ isomorphic to the action on the set $\Lambda$ of $k$-subsets of $\{1,\dots, m\}$.

If $G\cong\Sym(m)$, then Lemma~\ref{l: product fixed} implies that a non-trivial element of $G$ must move at least $2m-4$ points of $\Lambda$. We know that there exist non-trivial elements that move at most $6\lambda-6$ elements, and so we conclude that $m\leq 3\lambda-1$.

If $G\cong \Alt(m)$ with $m>5$, then Lemma~\ref{l: product fixed} implies that a non-trivial element of $G$ must move at least $3m-6$ points of $\Lambda$, and the same argument implies that $m\leq 2\lambda$.

If $G\cong \Alt(5)$, then Lemma~\ref{l: product fixed} implies that a non-trivial element of $G$ must fix at least $8$ points of $\Lambda$. We conclude that $6\lambda-6\geq 8$ and so $\lambda\geq 3$, and we are done.
\end{proof}

\section{Theorem C}\label{s: l3}

Our aim in this section is to classify puzzle groups arising from supersimple $2-(n,4,3)$ designs. Note, first, that such designs only occur for $n\equiv 0,1\pmod 4$ and $n\geq 8$. 

Throughout this section we let $\De$ be a supersimple $2-(n,4,3)$ design and set $G=\pi_{\infty}(\De)$. Since $\lambda=3$, we observe that all elementary moves are even permutations and so, by Lemma~\ref{Properties of D}, $G$ is a subgroup of $\Alt(n-1)$. We start by applying Theorem~A to this situation in which case we obtain the following lemma.

\begin{lem}\label{l: start 3}
\hspace{2em}
\begin{enumerate}
 \item if $n>11$, then $G$ is transitive;
 \item if $n>21$, then $G$ is primitive;
 \item if $n>50$, then $G=\Alt(n-1)$.
\end{enumerate}
\end{lem}


\subsection{Small \texorpdfstring{$n$}{n}}

The $2-(8,4,3)$ and $2-(9,4,3)$ designs are listed explicitly in \cite{Col}. Direct calculation then yields the following result.

\begin{lem}\label{l: small n}
The following statements holds:
\begin{enumerate}
\item There is a unique supersimple $2-(8,4,3)$ design, and its hole stabilizer is trivial.
\item There is a unique supersimple $2-(9,4,3)$ design, and its hole stabilizer, $G$, is transitive and imprimitive, with $G\cong \Alt(4)\wr C_2$.
\end{enumerate}
\end{lem}
\begin{proof}
Using the list in \cite{Col}, for $n=8$, we can see that there exists exactly one $2-(8,4,3)$ supersimple design. This designs is 
(isomorphic to) the Boolean quadruple system of order $8$, and so $\pi_{\infty}(\mathcal{D})$ is trivial.

For $n=9$ we can also check that exists exactly one $2-(9,4,3)$ supersimple design. A calculation using \cite{GAP} shows that
$\pi_{\infty}(\mathcal{D})\cong \Alt(4)\wr C_2$.
\end{proof}

Let us be explicit for the case $n=9$: it turns out that 
the only supersimple $2-(9,4,3)$ design is 
$$
\begin{array}{rcl}
\mathcal{D} & = & \{(1,2,3,4),(1,2,5,6),(1,2,7,8),(1,3,5,9),(1,3,6,7),(1,4,5,8),(1,4,7,9),\\[1.3ex]
			&   & (1,6,8,9),(2,3,5,7),(2,3,8,9),(2,4,5,9),(2,4,6,8),(2,6,7,9),(3,4,6,9),\\[1.3ex]
			&	& (3,4,7,8),(3,5,6,8),(4,5,6,7),(5,7,8,9)\}.\\[1.3ex]
\end{array}
$$

Next, a computer calculation of Professor Brendan McKay confirms that there are 28,893 supersimple $2-(12,4,3)$ designs; more computer calculations with \cite{GAP} confirm that all of these designs have hole stabilizer isomorphic to $\Alt(11)$, thus we assume that $13\leq n \leq 29$ from here on.

From here on we assume that $n\geq 13$. Lemma~\ref{l: start 3} implies, then, that $G$ is transitive.

\subsection{The imprimitive case}

Suppose that $G$ is transitive and preserves a system of imprimitivity with $\ell$ blocks of size $k$ (so $n-1=k\ell$). Lemma~\ref{l: start 3} implies that $n\leq 21$.  We know that $n-1$ cannot be prime, so this implies that $n\in\{13,16,17,21\}$.

\subsection{The primitive case}

In this section we assume that $G$ is primitive and not isomorphic to $\Alt(n-1)$. We know already, thanks to Lemma~\ref{l: start 3}, that $n\leq 50$ and, thanks to Lemma~\ref{l: small n}, that $n\geq 13$. We start by improving this.

\begin{lem}\label{L=3}
Suppose that $G$ is primitive. Then either $G\cong \Alt(n-1)$ or one of the following statements holds:
\begin{itemize}
\item $n=13$ and $G\in \{M_{12},M_{11}\}$;
\item $n=16$ and $G\in \{\mathrm{SL}_4(2), \Sym(6), \Alt(7), \Alt(6)\}$;
\item $n=17$ and $G$ is isomorphic to one of $18$ primitive groups in $2^4.\mathrm{SL}_4(2)$;
\item $n=28$ and $G\in \{\mathrm{PSp}_4(3)\rtimes C_2\}$;
\item $n=29$ and $G\in \{\mathrm{Sp}_6(2),\Sym(8)\}$.
\end{itemize}
\end{lem}
\begin{proof}
We know, by Theorem~\ref{Trivial E}, that there exist points $a,b \in \Omega$ such that $g=[\infty, a,b,\infty]$ is non-trivial and $a,b$ are collinear with $\infty$. Then $g$ is an element with support of size at most $6\lambda-6=12$. 

Now the list above contains all but one of the primitive groups on $n-1$ points which
\begin{enumerate}
 \item satisfy $13\leq n \leq 50$ with $n\equiv 0,1\pmod 4$;
 \item contain a non-trivial element with support at most $12$;
 \item are subgroups of $\Alt(n-1)$.
\end{enumerate}

Let us consider the missing entry which occurs when $n=13$ and $G\cong \PSL_2(11)$. It is easy to check that there are no non-trivial elements that fix more than 2 points. But now, by Theorem~\ref{Trivial E}, we can assume that there exists $g:=[\infty, a,b,\infty]$ which is not trivial and for which there exists $c$ such that $\{\infty, a,b,c\}\in\B$. But now observe that $g$ fixes $a$, $b$ and $c$, and so we have a contradiction.
\end{proof}

\begin{lem}\label{lem:5.4}
 $n\leq 17$.
\end{lem}
\begin{proof}
Suppose that $n>17$. Then $n\in \{28,29\}$. The three possible permutation groups 
given in Theorem~\ref{L=3} have precisely one non-trivial conjugacy class of 
elements of support at most $12$. In every case it is a conjugacy class of 
involutions with support exactly $12$. 

Using \cite{GAP} one can verify that if $G$ is one of these three permutation 
groups, $g, h\in G$ are two involutions of support $12$ and $\tau$ is one of 
the six disjoint transpositions whose product is $g$, then $\tau$ is not one 
of the six disjoint transpositions whose product is $h$.

Let $a\in \Omega\setminus\{\infty\}$ and consider the three lines connecting 
$\infty$ to $a$:
\[
 \{\infty, a, b_1, c_1\} \,\,  \{\infty, a, b_2, c_2\} \,\, \{\infty, a, b_3, c_3\}.
\]
Consider the permutation $[\infty, a, b_1, \infty]$. This is either trivial, 
or has support $12$. Note that in the latter case, this implies that the sets $\overline{\infty, a}$, $\overline{a,b_1}$ and $\overline{b_1,\infty}$ must overlap only in the set $\{\infty, a, b_1, c_1\}.$

On the other hand if $[\infty, a,b_1,\infty]$ is trivial, then one can check 
that one must have $b_2\in\overline{a, b_1}$ and one obtains that 
$[\infty, a,b_2,\infty]$ is trivial, likewise $[\infty, a,b_3,\infty]$.

By running the same argument starting with $b_2$ and $b_3$ in place of $b_1$, 
one sees that $[\infty, a,b_1,\infty]$ is trivial if and only if 
$[\infty, a,b_2,\infty]$ is trivial if and only if $[\infty, a, b_3,\infty]$ 
is trivial.

Now Theorem~\ref{Trivial E} implies that we can choose $a$ so that 
$[\infty, a,b_1,\infty]$ is not trivial. Thus the same is true of 
$[\infty, a,b_2,\infty]$. But now, note that both $[\infty, a, b_1,\infty]$ 
and $[\infty, a,b_2,\infty]$ include the transposition $(b_3,c_3)$. On 
the other hand $[\infty, a,b_1,\infty]$ includes the transposition 
$(b_2,c_2)$ which $[\infty, a, b_2, \infty]$ does not. This is a 
contradiction and we are done.
\end{proof}

Lemma~\ref{lem:5.4} completes the proof of Theorem~C. The remaining couple of results rule out some of the open possibilities from Theorem~C. The first of these results generalizes the idea of Lemma \ref{lem:5.4}.

\begin{lem}\label{lema 6(lambda - 1)}
Let $G:=\pi_{\infty}(\mathcal{D})$ a puzzle group, where 
$\mathcal{D}:=(\Omega, \mathcal{B})$ is a $2-(n, 4, \lambda)$ supersimple design.
Then, one of the next statements is true:
\begin{enumerate}
\item there exists a non trivial element of support strictly less than $6(\lambda - 1)$.
\item there exist two different elements $g,h$, both with cycle type $2^{3(\lambda - 1 )}$ and so that in their cycle decomposition they have a common transposition.
\end{enumerate}
\end{lem}

\begin{proof}
Suppose that (1) is not true, i.e.\ suppose that the unique element of $G$ with support less than $6 (\lambda-1)$ is 
the identity. Let $a \in \Omega \smallsetminus \{\infty\}$, and consider the $\lambda$ lines connecting 
$\infty$ to $a$:

\[  
\ell_1 := \{\infty, a, b_1, c_1 \}, \quad \ell_2 := \{\infty, a, b_2, c_2 \}, \quad  \ell_3 := \{\infty, a, b_3, c_3 \}, \quad \ldots ,
\quad \ell_{\lambda} := \{\infty, a, b_{\lambda}, c_{\lambda} \}. 
\]

Consider the permutations $g_r := [\infty, a , b_1 , \infty]$ for $r=1,\dots, \lambda$. Suppose first that $g_1$ is trivial. 
Observe that $(b_r,c_r)$ is a transposition in $[\infty, a]$,
and, according to supersimplicity $b_r \in \ell_{r}$, and only in $\ell_{r}$. So, $b_r \in \overline{a,b_1}$
or $b_r \in \overline{b_1, \infty}$. Then, $g_r := [\infty, a, b_r , \infty]$ is trivial, also.  Changing $b_1$ with $b_r$ and $b_r$ with $b_1$ we obtain that $g_r$ is trivial if 
and only if $g_1$ is trivial. 

Now, according to Theorem \ref{Trivial E} we can choose $a \in \Omega \smallsetminus \{\infty\}$
so that $g_1$ is not trivial.  Thus the same is true for $g_r$, for $r=2, 3, \ldots , \lambda$. Note, moreover, that the support of $g_i$ has at most $6(\lambda-1)$ elements, for $r=1,\dots, \lambda$. Then, since (1) is not true, we conclude that $|\supp(g_r)|= 6 (\lambda - 1)$ for $r = 2, \ldots , \lambda$. It follows that the elements $g_r$ have cycle type $2^{3(\lambda - 1 )}$. Now, note that $g_1$ and $g_2$ include the transposition $(b_3,c_3)$, and the elements $b_3, c_3$ can't 
appear in any other transposition of $g_1$ and $g_2$, so $g_1 \neq g_2$, and this proves $(2)$.
\end{proof}

\begin{lem}\label{primit:lambda=3}
Let $G:=\pi_{\infty}(\mathcal{D})$ a puzzle group, where 
$\mathcal{D}:=(\Omega, \mathcal{B})$ is a $2-(17, 4, 3)$ supersimple design.
If $G$ is primitive, then $G$ is not isomorphic to
any of the groups $ 2^4:D(2*5)$, $(\Alt(4) \times \Alt(4)):2$, $(2^4:5).4$, $\AGL_1(16):2$,
$\AGammaL_1(16)$. 
\end{lem}

\begin{proof}
A calculation with \cite{GAP} shows that all non-trivial elements in the listed groups have support at least $6(\lambda - 1 ) = 12$. Furthermore, these elements have cycle structure $2^{6}$. Also, for every pair of different elements $g,h$ of this cycle type, we have that every transposition in $g$ is different to every transposition in $h$. 

Now, applying Lemma \ref{lema 6(lambda - 1)} we get a contradiction.

\end{proof}

Combining this result with the earlier restrictions given in Lemma~\ref{L=3}, we find that there are precisely 14 possible hole-stabilizers when $n=17$. In the {\tt GAP} library, they are {\tt PrimitiveGroup(16,i)} for $i\in \{7,8,10,11,12,13,14,15,16,17,18,19,20,21\}$. The case $i=21$ corresponds to $\Alt(16)$.

More generally, for $n\leq 17$, the only known example when $G$ is primitive and not equal to $\Alt(n-1)$ occurs when $n=16$, and $G\cong \Sym(6)={\rm Sp}_4(2)$. This example was first noted in \cite{GGNS}, and then generalized to an infinite family in \cite{GGS}.

\subsection{The case when \texorpdfstring{$\lambda=4$}{Lambda=4}}

It would be interesting to see if it might be possible to extend the classification to include puzzle groups arising from supersimple $2-(n,4,4)$ designs. Note, first, that such designs only occur for $n\equiv 1\pmod 3$ and $n\geq 10$. 

If $n=10$, an easy counting argument confirms that a supersimple $2-(10,4,4)$ design is also a $3-(10,4,1)$ design. Now a computer calculation by Professor Brendan McKay confirms that there is only one such design, and its hole stabilizer is $\Sym(9)$. For the record, this design has 30 lines as follows:
\[
\begin{array}{l}
{[}5, 7 ,8, 9], [4, 6 ,8 ,9],[4 ,5 ,6, 7,], [3 ,6 ,7, 9],[3, 4, 5, 8],[2, 6, 7, 8], [2, 4, 5, 9], [2 ,3 ,8, 9],  \\ {[}2, 3, 5, 6], [2, 3, 4, 7], [1, 5, 6, 8], [1, 4, 7, 9], [1, 3, 7, 8], [1, 3, 5, 9], [1, 3 ,4 ,6], [1, 2, 6, 9], \\ {[}1, 2, 5, 7],[1, 2, 4, 8], [10, 5, 6, 9], [10, 4, 7, 8], [10, 3, 6, 8], [10, 3, 5, 7], [10, 3, 4, 9],\\ {[}10, 2, 7, 9],[10, 2, 5, 8],[10, 2, 4, 6],[10, 1, 8, 9],[10, 1, 6, 7],[10, 1, 4, 5], [10, 1, 2, 3].
\end{array}
\]

\end{document}